# Dissipative stabilization of linear input delay systems via dynamical state feedback controllers: an optimization based approach ⋆

**Qian Feng** * **Bo Wei** **

* *Department of Control and Computer Engineering, North China Electric Power University, Beijing, China (Corresponding Author emails: 52702714@ncepu.edu.cn, qfen204@aucklanduni.ac.nz)*
** *Department of Control and Computer Engineering, North China Electric Power University, Beijing, (e-mail: bowei@ncepu.edu.cn)*

**Abstract:** In this note, we present an effective solution to the stabilization of linear input delay systems subject to dissipative constraints while all the effect of input delay is compensated by a controller with novel structure. The method is inspired by the recent development in the mathematical treatment of distributed delays and predictor controllers, which are critical for the derivation of the solution. An important conceptual innovation is the use of a parameterized dynamical state feedback controller (DSFC), where the dimension of the controller equals the dimension of the control input. A sufficient condition for the existence of a dissipative DSFC is obtained via the Krasovskii functional approach, where the condition includes a bilinear matrix inequality (BMI). To solve the BMI, we apply an inner convex approximation algorithm which can be initialized based on an explicit construction of a predictor controller gain. The proposed DSFC can be considered as an extension of the classical predictor controller, thereby capable of compensating all the effects of the pointwise input delay while satisfying dissipative constraints. A numerical example is given to illustrate the effectiveness of our proposed methodology.

*Keywords:* Dissipativity; Dynamical State Feedback; Krasovskii Functional; Predictor Controller.

## 1. INTRODUCTION

It is well known Olbrot (1978); Krstic (2008) that a system $\dot{\boldsymbol{x}}(t) = A\boldsymbol{x}(t) + B\boldsymbol{u}(t-r)$ is always stabilizable by a predictor controller

$$\boldsymbol{u}(t) = K\,\mathrm{e}^{Ar}\,\boldsymbol{x}(t) + K\int_{-r}^{0} \mathrm{e}^{-A\tau}\,B\boldsymbol{u}(t+\tau)\mathsf{d}\tau \qquad (1)$$

for arbitrary $r > 0$ if $(A, B)$ is stabilizable in an LTI system sense. By using a predictor controller, the problem of stabilizing $\dot{\boldsymbol{x}}(t) = A\boldsymbol{x}(t) + B\boldsymbol{u}(t-r)$ is converted to a standard LTI state feedback problem, thereby eliminating the effect of the input delay. For existing works related to the design of predictor controllers for pointwise input delays, see Mondié and Michiels (2003); Fiagbedzi (2008); Tsubakino et al. (2016) and the references therein.

Though the use of (1) can compensate all the effects of input delay, it is not easy to add performance objectives to the design of (1) as the procedure is constructive. If we see (1) from the perspective of delay reduction Artstein (1982), the performance guaranteed by the equivalent LTI delay-free system does not translate to the performance of the resulting delay system with (1). This makes the optimal control with predictor controllers very challenging.

A potential solution to the aforementioned problem is that one may analyze the distributed delay (DD) in (1) directly in an optimization setting. For both time and frequency domain approaches, the predictor controller in (1) is an integral delay equation whose numerical implementation can lead to potential instability problems due to the hidden neutral term Mondié and Michiels (2003). The problems have been resolved in Kharitonov (2015) by the use of a dynamical state feedback controller (DSFC)

$$\dot{\boldsymbol{u}}(t) = (KB + X)\,\boldsymbol{u}(t) + (KA - XK)\left(\mathrm{e}^{Ar}\,\boldsymbol{x}(t) + \int_{-r}^{0} \mathrm{e}^{-A\tau}\,B\boldsymbol{u}(t+\tau)\mathsf{d}\tau\right) \qquad (2)$$

with a Hurwitz matrix $X$, which can totally compensate the input delay in $\boldsymbol{u}(t-r)$ and ensure safe numerical implementation. Since (2) is always constructible from (1), one may utilize the controller structure in (2) as the foundation for the stabilization of input delay systems with performance objectives, where $\boldsymbol{u}(t)$ will become a new state of the final closed-loop system (CLS).

There are two potential issues with using (2) for the stabilization of input delay systems with performance objectives. Firstly, one has to deal with the DD in (2) directly and treat the CLS with (2) as a linear system with a DD.

⋆ This work was partially supported by the NSFC (National Natural Science Foundation of China) grant 61903140; and by Grant 2020MS019 of Fundamental Research Funds for Central Universities, China.

Secondly, the structure in (2) is rather restricted by its relation to the predictor controller in (1). Hence a more general form of parametrization is preferable for a DSFC to ensure better performance. Thanks to the mathematical treatment on DDs we have developed in Feng and Nguang (2016), one can directly handle the delay kernels in $e^{-A\tau}B$ and use more functions for the construction of Krasovskii functionals (KFs). On the other hand, it is more straightforward to design a DSFC with general structure and performance objectives via the KF approach than most frequency domain approaches Michiels and Niculescu (2014); Iftar (2018); Apkarian and Noll (2018). Consequently, the idea in Feng and Nguang (2016) paves the way to develop a unified optimization framework to tackle input delay stabilization problems subject to dissipativity constraints or uncertainties.

In this paper, we propose an effective method for the stabilization of a linear input delay system with dissipative constraints via the Krasovskii functional approach. The stabilization is achieved via a novel DSFC which includes both pointwise and distributed delays. The controller generalizes the structure in (1) with a greater degree of freedom while being able to compensate the effects of input delay and simultaneously dealing with disturbance. To ensure the delay in $\boldsymbol{u}(t-r)$ is always fully compensated, the integral kernels of the KF always includes the kernels in $e^{-A\tau}B$ to ensure that the proposed controller always cover the case in (2). Sufficient conditions for the existence of controllers are obtained via the construction of a complete-type KF, where the conditions are expressed via matrix inequalities. To solve the resulting bilinear matrix inequality in the synthesis condition, an iterative algorithm is proposed based on the idea of inner convex approximation outlined in Dinh et al. (2012), where each iteration is expressed in terms of linear matrix inequalities (LMIs). The algorithm requires a feasible solution to start with, which can be provided based on the construction of (1) due to the connections between the proposed controller and (2). As a result, the initialization of the iterative algorithm does not require a separated theorem where the BMI is convexified and often with induced conservatism Feng and Nguang (2016). It is a unique advantage of the proposed method due to the use of a DSFC.

The paper is organized as follows: The synthesis problem is first formulated in Section 2 where the derivation of the CLS and the proposed DSFC are presented. The main results concerning dissipative synthesis are set out in Sections 3, which includes a theorem and an algorithm. Finally, a numerical example is tested in Section 4 prior to the conclusion.

**Notation :** The notations in this paper are presented as follows; $\mathbb{S}^n := \{X \in \mathbb{R}^{n\times n} : X = X^\top\}$; $\mathbb{R}^{n\times n}_{[n]} := \{X \in \mathbb{R}^{n\times n} : \mathrm{rank}(X) = n\}$; $\mathbb{C}(\mathcal{X};\mathbb{R}^n)$ with $\sup_{\tau\in\mathcal{X}} \|\boldsymbol{f}(\tau)\|_2$ is the Banach space of continuous functions with an uniform norm. $\mathbf{Sy}(X) := X + X^\top$ is the sum of a matrix with its transpose. A column vector containing a list of objects is denoted by $\mathbf{Col}_{i=1}^n x_i := \left[\mathbf{row}_{i=1}^n x_i^\top\right]^\top = \left[x_1^\top \cdots x_i^\top \cdots x_n^\top\right]^\top$. $*$ is applied to denote $[*]YX = X^\top YX$ or $X^\top Y[*] = X^\top YX$. The diagonal sum of two matrices are denoted by $X \oplus Y = \mathsf{Diag}(X,Y)$, respectively. Moreover, we use $\widetilde{\forall} x \in \mathcal{X}$ to denote the meaning for almost all $x \in \mathcal{X}$ with respect to the Lebesgue measure. Denote the Kronecker product by $\otimes$. Finally, we assume the order of matrix operations to be *matrix (scalars) multiplications* $> \otimes > \oplus > +$.

## 2. PROBLEM FORMULATION

The paper concerns the dissipative stabilization of

$$\begin{aligned}\dot{\boldsymbol{x}}(t) &= A\boldsymbol{x}(t) + B\boldsymbol{u}(t-r) + D_1\boldsymbol{w}(t), \quad \widetilde{\forall} t \geq t_0 \in \mathbb{R}\\ \boldsymbol{z}(t) &= C_1\boldsymbol{\chi}(t) + C_2\boldsymbol{\chi}(t-r)\\ &\quad + \int_{-r}^0 \widetilde{C}_3(\tau)\boldsymbol{\chi}(t+\tau)\mathrm{d}\tau + D_3\boldsymbol{w}(t), \quad r > 0\end{aligned} \tag{3}$$

with the supply rate function

$$\begin{gathered}\mathsf{s}(\boldsymbol{z}(t),\boldsymbol{w}(t)) = \begin{bmatrix}\boldsymbol{z}(t)\\ \boldsymbol{w}(t)\end{bmatrix}^\top \begin{bmatrix}\widetilde{J}^\top J_1^{-1}\widetilde{J} & J_2\\ * & J_3\end{bmatrix}\begin{bmatrix}\boldsymbol{z}(t)\\ \boldsymbol{w}(t)\end{bmatrix},\\ \widetilde{J}^\top J_1^{-1}\widetilde{J} \preceq 0, \;\; J_1^{-1} \prec 0, \;\; \widetilde{J} \in \mathbb{R}^{m\times m},\\ J_2 \in \mathbb{R}^{m\times q}, \;\; J_3 \in \mathbb{S}^q\end{gathered} \tag{4}$$

where $\boldsymbol{x}(t) \in \mathbb{R}^n$ is the state variable, $\boldsymbol{u}(t) \in \mathbb{R}^p$ stands for system's input with $\boldsymbol{\chi}(t) := \mathbf{col}\left(\boldsymbol{x}(t),\boldsymbol{u}(t)\right)$ and $\boldsymbol{z}(t) \in \mathbb{R}^m$ indicates the regulated output and a deterministic disturbance signal $\boldsymbol{w}(\cdot) \in \mathbb{L}^2\left([t_0,+\infty];\mathbb{R}^q\right)$. The dimensions of the state space matrices $A \in \mathbb{R}^{n\times n}$, $B \in \mathbb{R}^{n\times p}$, $D_1 \in \mathbb{R}^{n\times q}$ and $C_1, C_2, \widetilde{C}_3(\tau) \in \mathbb{R}^{m\times\nu}$, $D_3 \in \mathbb{R}^{m\times q}$ are determined by $m,n,p,q \in \mathbb{N}$ with $\nu := n+p$. Finally, we assume that $\exists K \in \mathbb{R}^{p\times n}$ such that $A + BK$ is Hurwitz.

*Remark 1.* The differential equation in (3) holds for almost all $t \geq t_0$ in the sense of Lebesgue measure since $\boldsymbol{w}(\cdot) \in \mathbb{L}^2\left([t_0,+\infty];\mathbb{R}^q\right)$ where $\mathbb{L}^2\left([t_0,+\infty];\mathbb{R}^q\right)$ follows the standard definition of $\mathbb{L}^p$ space.

*Remark 2.* The supply rate function in (3) is grounded in the quadratic constraints in Scherer et al. (1997), which features numerous performance criteria such as

- $\mathbb{L}^2$ gain performance: $J_1 = -\gamma I_m$, $\widetilde{J} = I_m$, $J_2 = \mathsf{O}_{m,q}$, $J_3 = \gamma I_q$ with $\gamma > 0$
- Strict Passivity: $J_1 \prec 0$, $\widetilde{J} = \mathsf{O}_m$, $J_2 = I_m$, $J_3 = \mathsf{O}_m$
- Sector constraints: $J_1 = \widetilde{J} = -I_m$, $J_2 = -\frac{1}{2}(\alpha+\beta)$, $J_3 = \alpha\beta I_m$ with $m = q$.

### 2.1 An example of DSFC

Using the classical predictor controller

$$\boldsymbol{u}(t) = K\boldsymbol{x}(t+r) = K\left(e^{Ar}\boldsymbol{x}(t) + \int_{-r}^0 e^{-A\tau}B\boldsymbol{u}(t+\tau)\mathrm{d}\tau\right) \tag{5}$$

can exponentially stabilize the nominal system in (3) with $\boldsymbol{w}(t) \equiv 0$, given $A + BK$ is Hurwitz for some $K \in \mathbb{R}^{p\times n}$. To secure stable numerical implementation Mondié and Michiels (2003) of the distributed delay in the controller, the DSFC in eq.(2) (Kharitonov, 2015, See eq.(10) for more information) is always constructible [1] where $X \in \mathbb{R}^{p\times p}$ is any Hurwitz matrix. The spectrum of the nominal system in (3) with (2) is

$$\Big\{s \in \mathbb{C} : \det\left(sI_n - A - BK\right)\det\left(sI_p - X\right) = 0\Big\} \tag{6}$$

which implies that exponential stability is always achievable for some $K \in \mathbb{R}^{p\times n}$ given $X \in \mathbb{R}^{p\times p}$ is Hurwitz.

[1] The controller structure in (Kharitonov, 2015, eq.(10)) is a generalization of the controller applied here

### 2.2 Formulation of Synthesis Problem

To handle DDs with an optimization framework, we need the following assumption which is similar to the Assumption 1 in Feng and Nguang (2016).

*Assumption 3.* There exist $C_3 \in \mathbb{R}^{m\times d\nu}$, $\Gamma \in \mathbb{R}^{n\times dp}$ and $\boldsymbol{f}(\cdot) \in \mathbb{C}^1\left([-r,0];\mathbb{R}^d\right)$ such that $\widetilde{C}_3(\tau) = C_3\left(\sqrt{\mathsf{F}}\boldsymbol{f}(\tau)\otimes I_\nu\right)$ and $\left[\mathsf{O}_{p\times n}, (KA - XK)\,\mathrm{e}^{-A\tau} B\right] = \Gamma\left(\sqrt{\mathsf{F}}\boldsymbol{f}(\tau)\otimes I_\nu\right)$, where $\boldsymbol{f}(\cdot)$ satisfies

$$\begin{aligned}&\exists \Pi \in \mathbb{R}^{d\times d}: \frac{\mathsf{d}\boldsymbol{f}(\tau)}{\mathsf{d}\tau} = \Pi \boldsymbol{f}(\tau)\\ &\mathsf{F}^{-1} = \int_{-r}^{0} \boldsymbol{f}(\tau)\boldsymbol{f}^\top(\tau)\mathsf{d}\tau \succ 0,\end{aligned} \tag{7}$$

and $X \in \mathbb{R}^{n\times n}$ is any Hurwitz matrix with any $K \in \mathbb{R}^{p\times n}$ rendering $A + BK$ Hurwitz.

*Remark 4.* Note that $\mathsf{F}^{-1}$ is well defined based on Theorem 7.2.10 in Horn and Johnson (2012).

*Remark 5.* Assumption 1 does not impose structural restriction on the matrices $A$ and $B$, since the identity $\frac{\mathsf{d}\,\mathrm{e}^{-A\tau} B}{\mathsf{d}\tau} = -A\,\mathrm{e}^{-A\tau} B$ is naturally compatible with (7). Meanwhile, the elements in $\boldsymbol{f}(\cdot)$ need to cover all the items in $\mathrm{e}^{-A\tau} B$ and $\widetilde{C}_3(\tau)$. Moreover, it is important to point out that one can add an unlimited number of functions to $\boldsymbol{f}(\cdot)$ satisfying (7), since the column dimension of $\boldsymbol{f}(\cdot)$ is unbounded. Finally, we emphasize that the use of $\sqrt{\mathsf{F}}$ does not affect the existence of $C_3$ and $\Gamma$ because $\sqrt{\mathsf{F}}$ has full rank.

To stabilize (3) for arbitrary large $r > 0$, we apply a DSFC

$$\dot{\boldsymbol{u}}(t) = K_1\boldsymbol{\chi}(t) + K_2\boldsymbol{\chi}(t-r) + \int_{-r}^{0} K_3F(\tau)\boldsymbol{\chi}(t+\tau)\mathsf{d}\tau \tag{8}$$

where $F(\tau) := \left(\sqrt{\mathsf{F}}\boldsymbol{f}(\tau)\otimes I_\nu\right)$ with $\boldsymbol{f}(\tau)$ and $\mathsf{F}$ given in Assumption 3, and $K_1; K_2 \in \mathbb{R}^{p\times\nu}$ and $K_3 \in \mathbb{R}^{p\times d\nu}$ are unknowns to be computed. Meanwhile, we assume the dynamics of the controller operating in a disturbance existing environment is

$$\begin{aligned}\dot{\boldsymbol{u}}(t) = K_1\boldsymbol{\chi}(t) + K_2\boldsymbol{\chi}(t-r) + \int_{-r}^{0} K_3F(\tau)\boldsymbol{\chi}(t+\tau)\mathsf{d}\tau\\ + D_2\boldsymbol{w}(t),\end{aligned} \tag{9}$$

which is more realistic to consider by the synthesis analysis. The value of $D_2 \in \mathbb{R}^{p\times q}$ here can be arbitrary but known. Note that the controller we want to construct in (8) does not contain $D_2\boldsymbol{w}(t)$. Once the matrix parameters in (8) are obtained, we implement (8) instead of (9). This is because $D_2\boldsymbol{w}(t)$ is caused by the operating environment, which is not part of the controller.

*Remark 6.* The structure of the controller in (8) generalizes (2) which is capable of compensating all the effect of $\boldsymbol{u}(t-r)$ while being significantly more general than (2). Therefore, (8) has greater potential to increase performance than (2). Finally, the use of (9) for the synthesis analysis ensures that the resulting controller in (8) possesses the ability to handle disturbance introduced by its operating environment.

By using Assumption 3 and combining equations in (3) and (9), the dynamics of the resulting CLS operating in a disturbances existing environment is denoted as

$$\begin{aligned}\dot{\boldsymbol{\chi}}(t) &= \begin{bmatrix} A & \mathsf{O}_{n\times p}\\ & K_1\end{bmatrix}\boldsymbol{\chi}(t) + \begin{bmatrix}\mathsf{O}_n & B\\ & K_2\end{bmatrix}\boldsymbol{\chi}(t-r)\\ &\int_{-r}^{0}\begin{bmatrix}\mathsf{O}_{n\times d\nu}\\ K_3\end{bmatrix}F(\tau)\boldsymbol{\chi}(t+\tau)\mathsf{d}\tau + \begin{bmatrix}D_1\\ D_2\end{bmatrix}\boldsymbol{w}(t),\\ \boldsymbol{z}(t) &= C_1\boldsymbol{\chi}(t) + C_2\boldsymbol{\chi}(t-r)\\ &\quad + \int_{-r}^{0} C_3F(\tau)\boldsymbol{\chi}(t+\tau)\mathsf{d}\tau + D_3\boldsymbol{w}(t),\\ \boldsymbol{\chi}(t_0+\theta) &= \boldsymbol{\psi}(\theta),\ \ \theta\in[-r,0],\end{aligned} \tag{10}$$

where $\boldsymbol{\psi}(\cdot) \in \mathbb{C}([-r,0];\mathbb{R}^\nu)$ is the initial condition. By observing the equations in (10), we can reach the following conclusions.

*Conclusion 1.* The functional differential equation in (10) is of retarded type, which satisfies the properties outlined in (Kharitonov, 2015, Theorem 2). Namely, if the controller in (8) can exponentially stabilize the nominal system in (3) with $\boldsymbol{w}(t) \equiv 0$, then (8) can be implemented by the approximation rules of numerical integration without hidden instability problems.

*Conclusion 2.* The structure of (8) incorporates (2) as a special case. This means one can always find appropriate parameters for (8) to exponentially stabilize (3) with $\boldsymbol{w}(t) \equiv 0$. On the other hand, the stability of any open-loop system in (3) stabilized by (2) can be analyzed by the proposed methodology in this paper, though the stability condition is only sufficient.

## 3. MAIN RESULTS ON DISSIPATIVE CONTROLLER DESIGN

The following Krasovskii stability criterion is utilized to verify the stability of the trivial solution of (10).

*Lemma 7.* Let $\boldsymbol{w}(t) \equiv \mathbf{0}_q$ in (10) and $r > 0$ be given, then the trivial solution $\boldsymbol{\chi}(t) \equiv \mathbf{0}_\nu$ of (10) is uniformly asymptotically (exponentially) stable with any $\boldsymbol{\psi}(\cdot) \in \mathbb{C}([-r,0];\mathbb{R}^\nu)$ if there exist $\epsilon_1;\epsilon_2;\epsilon_3 > 0$ and a differentiable functional $\mathsf{v} : \mathbb{C}([-r,0];\mathbb{R}^\nu) \to \mathbb{R}$ with $\mathsf{v}(\mathsf{0}_\nu(\cdot)) = 0$ such that

$$\epsilon_1\|\boldsymbol{\psi}(0)\|_2^2 \le \mathsf{v}(\boldsymbol{\psi}(\cdot)) \le \epsilon_2\|\boldsymbol{\psi}(\cdot)\|_\infty^2 \tag{11}$$

$$\widetilde{\forall} t \ge t_0,\ \tfrac{\mathsf{d}}{\mathsf{d}t}\mathsf{v}(\boldsymbol{\chi}_t(\cdot)) \le -\epsilon_3\|\boldsymbol{\chi}(t)\|_2^2 \tag{12}$$

for any $\boldsymbol{\psi}(\cdot) \in \mathbb{C}([-r,0];\mathbb{R}^\nu)$ in (10), where $\|\boldsymbol{\psi}(\cdot)\|_\infty^2 := \sup_{-r_\nu\le\tau\le 0}\|\boldsymbol{\psi}(\tau)\|_2^2$. Furthermore, $\boldsymbol{\chi}_t(\cdot)$ in (12) is defined by $\forall t \ge t_0$, $\forall\theta \in [-r,0]$, $\boldsymbol{\chi}_t(\theta) = \boldsymbol{\chi}(t+\theta)$ in which $\boldsymbol{\chi} : [t_0 - r_\nu, \infty) \to \mathbb{R}^\nu$ satisfies (10) with $\boldsymbol{w}(t) \equiv \mathbf{0}_q$.

**Proof.** See Corollary 1 in Feng et al. (2020).

The following definition of dissipativity is grounded in the original framework outlined in Willems (1972).

*Definition 8.* The system in (10) with $\mathsf{s}(\boldsymbol{z}(t),\boldsymbol{w}(t))$ is said to be dissipative if there exists a differentiable functional $\mathsf{v} : \mathbb{C}([-r,0];\mathbb{R}^\nu) \to \mathbb{R}$ such that

$$\widetilde{\forall} t \ge t_0,\quad \tfrac{\mathsf{d}}{\mathsf{d}t}\mathsf{v}(\boldsymbol{\chi}_t(\cdot)) - \mathsf{s}(\boldsymbol{z}(t),\boldsymbol{w}(t)) \le 0 \tag{13}$$

with $t_0 \in \mathbb{R}$, $\boldsymbol{z}(t)$ and $\boldsymbol{w}(t)$ in (10). Moreover, $\boldsymbol{\chi}_t(\cdot)$ in (13) is defined by the equality $\forall t \ge t_0$, $\forall\theta \in [-r,0]$, $\boldsymbol{\chi}_t(\theta) = \boldsymbol{\chi}(t+\theta)$ with $\boldsymbol{\chi}(t)$ satisfying (10).

The main results of this paper are presented as follows

*Theorem 9.* Let all the parameters in Assumption 3 be given, then the CLS in (10) with the supply rate function in (4) is dissipative, and the trivial solution of (10) with $\boldsymbol{w}(t) \equiv \mathbf{0}_q$ is uniformly asymptotically (exponentially) stable if there exist $P \in \mathbb{S}^\nu$, $Q \in \mathbb{R}^{\nu \times d\nu}$, $R \in \mathbb{S}^{d\nu}$, $S; U \in \mathbb{S}^\nu$ and $[K_1\ K_2\ K_3] \in \mathbb{R}^{p\times(2\nu+d\nu)}$ such that the following conditions hold,

$$\begin{bmatrix} P & Q \\ * & R + I_d \otimes S \end{bmatrix} \succ 0, \quad S \succ 0, \quad U \succ 0, \tag{14}$$

$$\widehat{\mathbf{\Phi}} = \mathbf{\Phi} + \mathsf{Sy}\left[\mathbf{P}^\top(\mathbf{A} + \mathbf{B}\mathbf{K})\right] \prec 0, \tag{15}$$

where $\mathbf{P} := [P\ \mathsf{O}_\nu\ Q\ \mathsf{O}_{\nu\times q}\ \mathsf{O}_{\nu\times m}]$ and

$$\begin{aligned} \mathbf{A} &= \begin{bmatrix} A & \mathsf{O}_{n\times\nu} & B & \mathsf{O}_{n\times d\nu} & D_1 & \mathsf{O}_{n\times m} \\ \mathsf{O}_{p\times n} & \mathsf{O}_{p\times\nu} & \mathsf{O}_p & \mathsf{O}_{p\times d\nu} & D_2 & \mathsf{O}_{p\times m} \end{bmatrix}, \\ \mathbf{B} &= \begin{bmatrix} \mathsf{O}_{n\times p} \\ I_p \end{bmatrix}, \ \mathbf{K} := \left[K_1\ K_2\ K_3\ \mathsf{O}_{p\times(q+m)}\right], \end{aligned} \tag{16}$$

$$\begin{aligned} \mathbf{\Phi} := \mathsf{Sy}&\left(\begin{bmatrix} Q \\ \mathsf{O}_{\nu\times\nu d} \\ R \\ \mathsf{O}_{(q+m)\times\nu d} \end{bmatrix} \left[F(0)\ -F(-r)\ -\widehat{\Pi}\ \mathsf{O}_{\nu d\times(q+m)}\right]\right) \\ &+ (S + rU) \oplus [-S] \oplus [-I_d \otimes U] \oplus J_3 \oplus J_1^{-1} \\ &+ \mathsf{Sy}\left(\begin{bmatrix} \mathsf{O}_{(2\nu+\nu d)\times m} \\ J_2^\top \\ I_m \end{bmatrix} [\mathbf{\Sigma}\ \mathsf{O}_m]\right). \end{aligned} \tag{17}$$

with $\mathbf{\Sigma} = [C_1\ C_2\ C_3\ D_3]$ and $\widehat{\Pi} := \sqrt{\mathsf{F}}\Pi\sqrt{\mathsf{F}^{-1}} \otimes I_\nu$. The number of decision variables of Theorem 9 is $(0.5d^2 + d + 1.5)\nu^2 + (0.5d + pd + 2p + 1.5)\nu \in \mathcal{O}(d^2\nu^2)$ with $\nu = n + p$.

**Proof.** The proof is based on the construction of the KF

$$\begin{aligned} v(\boldsymbol{\chi}_t(\cdot)) := [*] \begin{bmatrix} P & Q \\ * & R \end{bmatrix} &\begin{bmatrix} \boldsymbol{\chi}(t) \\ \int_{-r}^0 F(\tau)\boldsymbol{\chi}(t+\tau)\mathsf{d}\tau \end{bmatrix} \\ + \textstyle\int_{-r}^0 \boldsymbol{\chi}^\top(t+\tau)\left[S + (\tau + r)U\right]&\boldsymbol{\chi}(t+\tau)\mathsf{d}\tau \end{aligned} \tag{18}$$

where $\mathbb{R}^{\nu\times\nu d} \ni F(\tau) = \sqrt{\mathsf{F}}\boldsymbol{f}(\tau) \otimes I_\nu$ and $P \in \mathbb{S}^\nu$, $Q \in \mathbb{R}^{\nu\times d\nu}$, $R \in \mathbb{S}^{d\nu}$, $S; U \in \mathbb{S}^\nu$.

The detail of the proof are omitted here due to limited space, and it will be presented in the journal version of this paper.

*Remark 10.* Theorem 9 can be extended to optimize the generalized $\mathcal{H}_2$ objective.[2] The corresponding synthesis constraints can be derived considering Appendix C.7 of Briat (2014) without difficulties.

*Remark 11.* The inequality in (15) is bilinear because of $\mathsf{Sy}\left(\mathbf{P}^\top\mathbf{B}\mathbf{K}\right)$. On the other hand, the BMI in (16) can be reformulated as an LMI subject to rank constraints which is still nonconvex.

The convexification techniques in Feng and Nguang (2016) based on Projection Lemma can be applied to $\mathsf{Sy}\left(\mathbf{P}^\top\mathbf{B}\mathbf{K}\right)$ at the expense of adding slack variables subject to potential conservatism. However, the problem we are solving here has a unique feature that

$$\begin{aligned} K_1 &= \left[(KA - XK)\,\mathrm{e}^{Ar}\ \ KB + X\right], \\ K_2 &= \mathsf{O}_{p\times\nu}, \quad K_3 = \Gamma \end{aligned} \tag{19}$$

[2] $D_3 = \mathsf{O}_{m\times q}$ is required to ensure a finite generalized $\mathcal{H}_2$ gain Briat (2014).

can always be considered as a feasible solution to the unknowns in (9) where $K, X$ and $\Gamma$ are given in Assumption 3. In consequence, this allows us to apply the inner convex approximation algorithm in Dinh et al. (2012) to solve the BMI in (16) whose initialization can be achieved via (19).

To apply the algorithm in Dinh et al. (2012), consider

$$\begin{aligned} \mathbb{S}^{\ell\times\ell} \ni \Delta\left(\mathbf{G}, \widetilde{\mathbf{G}}, \mathbf{N}, \widetilde{\mathbf{N}}\right) &:= [*]\left[Z \oplus (I_n - Z)\right]^{-1}\begin{bmatrix} \mathbf{G} - \widetilde{\mathbf{G}} \\ \mathbf{N} - \widetilde{\mathbf{N}} \end{bmatrix} \\ &+ \mathsf{Sy}\left(\widetilde{\mathbf{G}}^\top\mathbf{N} + \mathbf{G}^\top\widetilde{\mathbf{N}} - \widetilde{\mathbf{G}}^\top\widetilde{\mathbf{N}}\right) + \mathbf{T} \end{aligned} \tag{20}$$

with $Z \oplus (I_n - Z) \succ 0$ satisfying

$$\begin{aligned} \forall \mathbf{G}; \widetilde{\mathbf{G}} \in \mathbb{R}^{n\times\ell},\ \forall \mathbf{N}; \widetilde{\mathbf{N}} \in \mathbb{R}^{n\times\ell},\ \ &\mathbf{T} + \mathsf{Sy}\left(\mathbf{G}^\top\mathbf{N}\right) \\ = \Delta\left(\mathbf{G}, \mathbf{G}, \mathbf{N}, \mathbf{N}\right) &\preceq \Delta\left(\mathbf{G}, \widetilde{\mathbf{G}}, \mathbf{N}, \widetilde{\mathbf{N}}\right). \end{aligned} \tag{21}$$

This shows that $\Delta\left(\mathbf{G}, \widetilde{\mathbf{G}}, \mathbf{N}, \widetilde{\mathbf{N}}\right)$ is a psd-overestimate function Dinh et al. (2012) of $\acute{\Delta}(\mathbf{G}, \mathbf{N}) = \mathbf{T} + \mathsf{Sy}\left[\mathbf{G}^\top\mathbf{N}\right]$ with respect to the parameterization

$$\mathsf{Col}\left(\mathbf{vec}(\widetilde{\mathbf{G}}), \mathbf{vec}(\widetilde{\mathbf{N}})\right) = \mathsf{Col}\left(\mathbf{vec}(\mathbf{G}), \mathbf{vec}(\mathbf{N})\right). \tag{22}$$

Now let $\mathbf{T} = \mathbf{\Phi}$, $\mathbf{G} = \mathbf{P}$ and $\widetilde{P} \in \mathbb{S}^n$, $\widetilde{Q} \in \mathbb{R}^{n,dn}$

$$\begin{aligned} \widetilde{\mathbf{G}} = \widetilde{\mathbf{P}} &:= \left[\widetilde{P}\ \mathsf{O}_\nu\ \widetilde{Q}\ \mathsf{O}_{\nu\times q}\ \mathsf{O}_{\nu\times m}\right], \\ \mathbf{\Lambda} = [P\ Q]\,,\ \widetilde{\mathbf{\Lambda}} &:= \left[\widetilde{P}\ \widetilde{Q}\right],\ \mathbf{N} = \mathbf{B}\mathbf{K},\ \ \widetilde{\mathbf{N}} = \mathbf{B}\widetilde{\mathbf{K}} \\ \widetilde{\mathbf{K}} &= \left[\widetilde{K}_1\ \widetilde{K}_2\ \widetilde{K}_3\ \mathsf{O}_{p\times(q+m)}\right] \in \mathbb{R}^{p\times(2\nu+d\nu+q+m)} \end{aligned} \tag{23}$$

in (20) with $\ell =: 2\nu + d\nu + q + m$ and $Z \oplus (I_n - Z) \succ 0$ and $\mathbf{\Phi}$, $\mathbf{K}$ are given in Theorem 9.

Applying the matrix inequality in (21) with (23) to (16) produces

$$\begin{aligned} \widehat{\mathbf{\Phi}} &\preceq \mathbf{\Phi} + \mathsf{Sy}\left(\mathbf{P}^\top\mathbf{A} + \widetilde{\mathbf{P}}^\top\mathbf{N} + \mathbf{P}^\top\widetilde{\mathbf{N}} - \widetilde{\mathbf{P}}^\top\widetilde{\mathbf{N}}\right) \\ &+ \left[\mathbf{P}^\top - \widetilde{\mathbf{P}}^\top\ \ \mathbf{N}^\top - \widetilde{\mathbf{N}}^\top\right]\left[Z \oplus (I - Z)\right]^{-1}[*] \prec 0, \end{aligned} \tag{24}$$

where we assume that the upper bound of $\widehat{\mathbf{\Phi}}$ here is negative defined. Furthermore, $Z \in \mathbb{S}^\nu_{\succ 0}$ is a new decision variable. Applying Schur complement to (24) concludes that (24) holds if and only if (25) holds which now can be handled by standard SDP solvers provided that the values of $\widetilde{\mathbf{P}}$ and $\widetilde{\mathbf{K}}$ are given.

By compiling all the aforementioned procedures according to the expositions in Dinh et al. (2012), Algorithm 1 can be constructed as follows. To present a clear picture for the algorithm, we use vector $\mathbf{x}$ to contain all the decision variables in $R, S, U$. Furthermore, $\rho_1, \rho_2$ and $\varepsilon$ are given constants to achieve regularizations and determine error tolerance, respectively. To initialize the iterative algorithm, one can simply plug in the values in (19) where the parameters can be assigned via the construction of a predictor controller for $\dot{\boldsymbol{x}}(t) = A\boldsymbol{x}(t) + B\boldsymbol{u}(t - r)$.

## 4. NUMERICAL EXAMPLES

The numerical example is computed with Mosek (2022) programed with the optimization parser Yalmip Löfberg (2004) in Matlab.

Consider (3) with the parameters

$$
\begin{bmatrix} \widehat{\mathbf{\Phi}} + \mathsf{Sy}\left(\widetilde{\mathbf{P}}^\top \mathbf{A} + \widetilde{\mathbf{P}}^\top \mathbf{N} + \mathbf{P}^\top \widetilde{\mathbf{N}} - \widetilde{\mathbf{P}}^\top \widetilde{\mathbf{N}}\right) & \mathbf{P}^\top - \widetilde{\mathbf{P}}^\top & \mathbf{N} - \widetilde{\mathbf{N}}^\top \\ * & -Z & \mathsf{O}_n \\ * & * & Z - I_n \end{bmatrix} \prec 0 \tag{25}
$$

**Algorithm 1:** Inner Convex Approximation

**begin**
- **solve** Theorem 9 with (19) **return** $P$, $Q$
- **solve** Theorem 9 with $P$, $Q$ **return** $\mathbf{K}$.
- **update** $\widetilde{\mathbf{\Lambda}} \longleftarrow \mathbf{\Lambda}$, $\widetilde{\mathbf{K}} \longleftarrow \mathbf{K}$,
- **solve** $\min_{\mathbf{x},\mathbf{\Lambda},\mathbf{K}} \mathsf{tr}\left[\rho_1[*]\left(\mathbf{\Lambda} - \widetilde{\mathbf{\Lambda}}\right) + \rho_2[*]\left(\mathbf{K} - \widetilde{\mathbf{K}}\right)\right]$
  subject to the LMIs in (14), (23) and (25),
  **return** $\mathbf{\Lambda}$ and $\mathbf{K}$
- **while** $\dfrac{\|\mathbf{vec}(\mathbf{\Lambda}, \mathbf{K}) - \mathbf{vec}(\widetilde{\mathbf{\Lambda}}, \widetilde{\mathbf{K}})\|_\infty}{\|\mathbf{vec}(\widetilde{\mathbf{\Lambda}}, \widetilde{\mathbf{K}})\|_\infty + 1} \geq \varepsilon$ **do**
  - **update** $\widetilde{\mathbf{\Lambda}} \longleftarrow \mathbf{\Lambda}$, $\widetilde{\mathbf{K}} \longleftarrow \mathbf{K}$;
  - **solve** $\min_{\mathbf{x},\mathbf{\Lambda},\mathbf{K}} \mathsf{tr}\left[\rho_1[*]\left(\mathbf{\Lambda} - \widetilde{\mathbf{\Lambda}}\right) + \rho_2[*]\left(\mathbf{K} - \widetilde{\mathbf{K}}\right)\right]$
    subject to the LMIs in (14), (23) and (25),
    **return** $\mathbf{\Lambda}$ and $\mathbf{K}$;
- **end**

**end**

$$
\begin{aligned}
&A = \begin{bmatrix} -1 & 1 \\ 0 & 0.1 \end{bmatrix}, \; B = \begin{bmatrix} 0 \\ 1 \end{bmatrix}, \; D_1 = \begin{bmatrix} 0.1 \\ -0.1 \end{bmatrix} \\
&C_1 = \begin{bmatrix} -0.3 & 0.4 & 0.1 \\ -0.3 & 0.1 & -0.1 \end{bmatrix}, \; C_2 = \begin{bmatrix} 0 & 0.2 & 0 \\ -0.2 & 0.1 & 0 \end{bmatrix}, \\
&\widetilde{C}_3(\tau) = \begin{bmatrix} 0.2 + 0.1\,\mathrm{e}^{\tau} & 0.1 & 0.12\,\mathrm{e}^{3\tau} \\ -0.2 & 0.3 + 0.14\,\mathrm{e}^{2\tau} & 0.11\,\mathrm{e}^{3\tau} \end{bmatrix}, \\
&D_2 = 0.12, \;\; D_3 = \mathsf{Col}(0.14, 0.1).
\end{aligned} \tag{26}
$$

We first need to calculate an initial value for the controller gain $\mathbf{K}$ based on the construction of (2) and (5). Consider the fact that

$$
\mathrm{e}^{-A\tau} B = \mathsf{Col}\left[(10/11)\,\mathrm{e}^{-0.1\tau} - (10/11)\,\mathrm{e}^{\tau}, \; \mathrm{e}^{-0.1\tau}\right] \tag{27}
$$

and the functions inside of $\widetilde{C}_3(\tau)$, we apply

$$
\boldsymbol{f}(\tau) = \mathsf{Col}\left(1, \mathrm{e}^{\tau}, \mathrm{e}^{2\tau}, \mathrm{e}^{3\tau}, \mathrm{e}^{-0.1\tau}\right) \tag{28}
$$

associated with $\Pi = 0 \oplus 1 \oplus 2 \oplus 3 \oplus (-0.1)$ as the basis to denote all the distributed terms in (26)–(27).

Let $X = -0.1$ in (19) to apply Algorithm 1 as $X = -0.1$ is Hurwitz. Then one can obtain

$$
\begin{aligned}
C_3 = &\left[\begin{smallmatrix} 0.2 & 0.1 & 0 & 0.1 & 0 & 0 & 0 & 0 & 0 \\ -0.2 & 0.3 & 0 & 0 & 0 & 0 & 0 & 0.14 & 0 \end{smallmatrix}\right. \\
&\left.\begin{smallmatrix} 0 & 0 & 0.12 & 0 & 0 & 0 \\ 0 & 0 & 0.11 & 0 & 0 & 0 \end{smallmatrix}\right]\left(\sqrt{\mathsf{F}^{-1}} \otimes I_3\right) \\
K_1 = &\,[0.0235 \;\; -0.2629 \;\; -0.5173]\,, \; K_2 = \mathbf{0}_3^\top \\
K_3 = \Gamma = &\left[\mathbf{0}_5^\top \;\; -0.4295 \;\; \mathbf{0}_8^\top \;\; -0.1789\right]\left(\sqrt{\mathsf{F}^{-1}} \otimes I_3\right)
\end{aligned} \tag{29}
$$

by (19) and Assumption 3, which gives an initial value of $\mathbf{K}$ to apply Algorithm 1.

Applying Theorem 9 to (26) with (29) yields a feasible solution and shows the controller can achieve $\min \gamma = 0.49425$. Next, resubstitute the resulting $P$,$Q$ produced by the previous step into Theorem 9 which produces a new $\mathbf{K}$. Consequently, we obtain

$$
\begin{aligned}
K_1 &= [0.0294 \;\; -0.2732 \;\; -0.5098]\,, \\
K_2 &= [0.0057 \;\; -0.0053 \;\; -0.0077]\,, \\
K_3 &= \left[-0.0039 \; 0.0014 \; -0.3878 \; 0.0029 \cdots \right. \\
&\quad 0.0014 \; -0.2378 \; 0.0077 \; -0.006 \cdots \\
&\quad -0.1538 \; 0.0079 \; -0.0115 \; -0.0989 \cdots \\
&\quad \left. -0.0089 \; 0.0027 \; -0.4048\right]
\end{aligned} \tag{30}
$$

which guarantees $\gamma = 0.49227$. Now use (30) with the associated $P$ and $Q$ for $\widetilde{\mathbf{\Lambda}}$ and $\widetilde{\mathbf{K}}$ to use Algorithm 1. The results produced by Algorithm 1 with $\rho_1 = 0.01$, $\rho_2 = 0.01$ are summarized in Table 1 in which NoI standards for number of iterations. Note that here we prescribe a very small $\varepsilon = 10^{-10}$ to control the value of NoI. Furthermore, we do not present the resulting controller gains due to the limitation of space.

| $\min \gamma$ | 0.481 | 0.4714 | 0.46398 | 0.45749 |
|---|---|---|---|---|
| NoI | 100 | 200 | 300 | 400 |

Table 1. $\min \gamma$ produced by different iterations

Clearly, the results in Table 1 demonstrate that more iterations lead to better $\min \gamma$ value at the expense of more computational time. In addition, it also shows the advantage of using Algorithm 1 to calculate a controller gain with better performance.

To verify the validity of our controllers, we apply the spectrum method in Breda et al. (2015) to calculate the characteristic roots of all the resulting CLSs. It shows that the spectral abscissae of all the CLSs are negative, which implies that all the CLSs are exponentially stable. Finally, as all the resulting CLSs are of retarded type, the controller in (8) satisfies the properties elucidated in Conclusion 1, which means it can always be implemented numerically without causing instability problems.

## 5. CONCLUSION

We have presented an effective solution to the dissipative stabilization of the input delay system (3) by the use of a novel dynamical feedback controller (8). The key components of the proposed method consist of the decomposition in Assumption 3, the novel controller in (8), the construction of the KF in (18) and the inner convex approximation in Algorithm 1. The proposed method has attained our original goal to design a controller which can compensate all the effects of the input delay while guaranteeing performance objectives via the incorporation of the supply rate function in (4). Moreover, the design has also incorporated potential disturbance caused by the operating environment, hence the novel controller (8) is robust and able to fulfill performance requirements related to disturbance. As for Algorithm 1, its initial values can always be supplied by the gain of a predictor controller. On account of the unique structure in (10) where the state space parameters of (3) are not directly overlapped with the controller gains in (8), it is possible to design a resilient controller in the form of (8) for (3) with uncertainties. This is a crucial point for future research as the design of resilient controllers is quite challenging for uncertain systems.